%% file: Laurain-CMC-FB.tex
\documentclass[reqno]{amsart}

\usepackage[babel=true,kerning=true]{microtype}

\usepackage{amsthm,leqno}

\usepackage{amsfonts,amsmath,amssymb,enumerate}
\usepackage{graphicx}
\usepackage{tikz}
\usetikzlibrary{decorations.markings}
\usetikzlibrary{plotmarks}
\usetikzlibrary{patterns}

\newcommand{\cmc}{constant mean curvature}
\newcommand{\ds}{\displaystyle}

\newcommand{\be}{\begin{equation*}}
\newcommand{\ee}{\end{equation*}}
\newcommand{\beq}{\begin{equation}}
\newcommand{\eeq}{\end{equation}}
\newcommand{\begincal}{\begin{eqnarray*}}
\newcommand{\fincal}{\end{eqnarray*}}
\newcommand{\lbeq}{\begin{leqno}}
\newcommand{\leeq}{\end{leqno}}

\newtheorem{thm}{Theorem}[section]
\newtheorem{lemma}{Lemma}[section]
\newtheorem{cor}{Corollary}[section]

\newtheorem{defi}{Definition}[section]

\newtheorem{rem}{Remark}[section]

\newcommand{\eps}{{\varepsilon}}
\newcommand{\R}{{\mathbb R}}
\newcommand{\C}{{\mathbb C}}
\newcommand{\N}{{\mathbb N}}
\newcommand{\D}{{\mathbb D}}
\newcommand{\h}{{\mathbb H}}

\title{Asymptotic analysis for surfaces with large constant mean curvature and free boundaries}
\author{Paul Laurain}

\begin{document}

\begin{abstract} 
We prove that simply connected $H$-surfaces with bounded area  and free boundary in a domain necessarily concentrate at a critical point of the mean curvature of the boundary of this domain.
\end{abstract}
\maketitle

\section*{Introduction}

The aim of this article is to understand the asymptotic behaviour of sequences of surfaces with large constant mean curvature and free boundaries. These surfaces arise naturally in the partitioning problem which consists in dividing a domain into two parts of prescribed volumes by a surface of minimal area. The existence of solutions of this problem is given by the geometric measure theory (see for instance Morgan \cite{Morgan}). However we get no information about the topology of such surfaces, except in the case of strictly convex domains, where we know that such a surface is connected and we get some bounds  on the number of components of its boundary as well as its genus, see Ros \& Vergasta \cite{RosVergasta}. Moreover it is conjectured that, in this case, the surface is homeomorphic to a disk, see  Ritore \& Ros \cite{RitoreRos}.\\

In the following, we let $\Omega$ be a smooth domain of $\R^3$ and we will consider $H$-surface as a map $u\in C^2(\overline{\D},\R^3)$ where
\be
\D=\{ z\in\R^2 \hbox{ s.t. } \vert z \vert <1\} 
\ee
which is an immersion and which satisfies
\beq
\label{libre1}
 \begin{cases}
\Delta u=-2H\, u_x\wedge u_y, \\
\langle u_x,u_y\rangle =\vert u_x\vert-\vert u_y\vert =0,\\
u(z) \in \partial \Omega \hbox{ for all } z\in\partial\D ,\\
\partial_{\nu}u(z)\bot T_{u(z)}\partial \Omega \hbox{ for all } z\in\partial\D ,
\end{cases}
\eeq
where $\ds \Delta=-\frac{\partial^2 }{\partial x^2}-\frac{\partial^2 }{\partial y^2}$.\\

Then $u(\overline{D})$ is a regular surface of constant mean curvature $H$ with boundary contained in $\partial\Omega$ and which meets $\partial\Omega$ orthogonally.\\

The first result of existence of solutions of (\ref{libre1}) is due to Struwe \cite{Struwe88}, which finds solutions in domains diffeomorphic to a ball using a parabolic version of our equation. Another idea to find solutions of (\ref{libre1}) in a general domain  is to look for solutions with large mean curvature (i.e. with a small diameter). In fact in this case the topology of the domain play no role and the geometry is under control.\\

This intuition was confirmed by Fall \cite{Fall07}. However, the existence of such solutions is subject to a local condition on the  curvature of the boundary of $\Omega$. Fall proved in \cite{Fall07} the following : given any smooth domain $\Omega\subset \R^3 $ and $ p \in \partial \Omega $ be a non-degenerate critical point of the mean curvature of $ \partial \Omega $,There exists a family of solutions $ u^\eps \in C^2 (\overline{\D}, \R^3) $ of (\ref{libre1})  for $ H =  \frac{1}{\eps} $ such that $u^\eps (\D) $ is embedded and $ \Vert u^\eps -p   \Vert_\infty \rightarrow0 $ when $ \eps \rightarrow 0 $. Moreover $ \frac{1}{\eps} u^\eps $, correctly translated, converges to an hemisphere of radius $1$.\\

This result is similar to the result of Ye \cite{Ye91} concerning the existence of closed surfaces with constant mean curvature in a curved manifold. Indeed they are similar in their statement but also in the method of proof which takes a solution of the limit equation, here one hemisphere, and tries to perturb it via the implicit functions theorem. This remark done,  the question of the necessity of the condition that $p$ is a critical point of the mean curvature of $\partial\Omega$ comes naturally. A first answer in this direction is provided by Fall \cite{Fall10}. Indeed he shows that the solutions to the problem of partitioning, which are equivalent to the isoperimetric problem solutions in this context, converge to a point of maximal mean curvature when their volume tends to zero.\\

This theorem is similar to  the result of Druet \cite{Druet02} concerning the location of small isoperimetric domains in a curved manifold. Indeed Druet proved that these domains are near global maxima of the scalar curvature­. In \cite{laurain10}, we proved under suitable assumptions that surfaces of large constant mean curvature and small diameter in a $3$-dimensional manifold are necessarily located near a critical point of the scalar curvature. Here we show under reasonable assumptions  that surfaces of large constant mean curvature with boundary included in $\partial\Omega$ and meeting $\partial\Omega$ orthogonally are necessarily located near a critical point of the mean curvature of $\partial\Omega$. First, we assume that the diameter is controlled in order to avoid solutions that collapse along some geodesics (such examples were constructed by Mahmoudi and Fall \cite{FallMahmoudi}). Second, we assume that the area is controlled to avoid an infinity of bubbles. Then we prove the following theorem.
\begin{thm}
\label{Laurain10ter}
Let $\Omega$  be a smooth domain of $\R^3$ and a sequence of embedded surfaces  $\Sigma^\eps$ in $\Omega$ satisfying the following assumptions :
\begin{enumerate}[(i)]
\item  $\partial\Sigma^\eps \subset \partial\Omega$ and  $\Sigma^\eps$ and $\partial\Omega$ meets orthogonally,
\item  $\Sigma^\eps$ has \cmc\ equal to $\frac{1}{\eps}$,
\item the diameter and the area of $\Sigma^\eps$ are respectively a $O(\eps)$ and a $O(\eps^2)$.
\end{enumerate}
Then, up to a subsequence, $\Sigma^\eps$ converge to $p\in \partial\Omega$ which is a critical point of the mean curvature of $\partial\Omega$.
\end{thm}
This theorem can be explained in the following way: given $\Omega\subset\R^3$ a smooth domain, for any $\delta>0$, any $C>0$, there exists $\eps_0 >0$ such that any embedded surface  orthogonal to the boundary $\Sigma$ of constant mean curvature $\frac{1}{\eps}$ with $\eps<\eps_0$, $\text{diameter}(\Sigma)\leq C\eps$, $\text{Area}(\Sigma)\leq C\eps^2$, satisfies that $\Sigma\subset B(p,\delta)$ for some critical point $p\in \partial\Omega$ of the mean curvature of $\partial\Omega$.\\

Note that the bound on the diameter and the area are scale invariant with respect to the mean curvature.

This article is organized as follows. In the first section we remind some useful results about regularity of \cmc\ surfaces with free boundaries. In the second section we remind the classification of the solution of the constant mean curvature equation on the whole plane and we extend it to domain like disk or half-plane. Finally in a third section  we give a proof of the theorem, dividing it in three parts; first we perform a blow-up analysis decomposing our sequence in a sum of spheres and hemispheres; then we insure the existence of at least one hemisphere in the decomposition using notably the Aleksandrov reflexion principle, finally we achieve the proof applying the balancing formula. The main difficulty is to understand precisely the asymptotic behaviour of our sequence of surfaces $\Sigma^\eps$ on the boundary of $\Omega$. Some technical lemmas are postponed to the appendix.\\

\medskip {\bf Acknowledgements} : I would to thank Frank Pacard for submitted me this question and the work of Fall. It is also my pleasure to express my deep thanks to my thesis advisor O.Druet for his encouragement and helpful comments during the preparation of this work.
\section{Regularity and  {\it a priori}  estimates on \cmc\ surfaces with free boundaries }
\label{sregul}
In this section we give a general result on the regularity of constant mean curvature surfaces with free boundaries, the reader will find all the details in Chapter 7 of \cite{DierkesII}.
\begin{thm}
\label{regulb}
Let $\Omega$ be a $C^{m,\alpha}$ domain of $\R^3$ with $m\geq 3$ and $\alpha>0$,  then every solution of (\ref{libre1}) is  $C^{m,\alpha}$.
\end{thm}
The solutions inherit of the regularity of the $\Omega$ provided it is sufficiently smooth.\\

The proof of this result is divided into three steps. A first shows, using the isoperimetric inequality for surfaces, that the solutions are $C^{0,\eta}$ up to the boundary. Then, using  {\it a priori } estimates in the spaces $H ^{k, p}$, we deduce the $C^{1,\frac{1}{2}}$ regularity up to the boundary. Finally, using a classical argument of bootstrap, we obtain that the solutions are smooth inside and inherit of the regularity of the domain up to the boundary as soon as it is at least $ C ^{3, \alpha} $. \\

We give here the {\it a priori} estimate which is the keystone of the second step and that will be used later.
\begin{thm} Let $\Omega$ be a smooth domain, whose metric of the boundary  will be denoted by $g$, and $u$  be a solution of (\ref{libre1}). We assume that $u$ belongs to $C^{0,\eta}(\overline{\D})$. Then, for every open set  $U$ of $\overline{\D}$ and every $2<p<+\infty $, there exists a constant $c$ depending only on  $\Vert g\Vert_3 $, $U$, $p$ ,$\int_U \vert \nabla u\vert^2 dz$ and the modulus of continuity of   $u$ such that 
$$\int_{U} \vert \nabla u \vert^p dz < c.$$ 
\end{thm}
This estimate and the standard elliptic theory lead to uniform bounds of the type
$$\Vert u \Vert_{2+\eta,U} < c,$$
where $c$ depends only on $\Vert g\Vert_3 $, $U$, $p$ ,$\int_U \vert \nabla u\vert^2 dz$ and the modulus of continuity of   $u$. 
\begin{rem}
\label{estimreg}
In particular, we note that from any sequence of solutions whose gradient is uniformly bounded on an open set $ U $ of $ \overline{\D} $, we can extract a subsequence which converges uniformly in $C^2(U)$.
\end{rem}
\section{Classification of solution of the limit equation}
We start by remind a crucial result of Brezis and Coron \cite{BC2} which states that the only solutions of
\be
\Delta u = -2\, u_x \wedge u_y \hbox{ on } \R^2
\ee 
with bounded energy are exactly, up to a conformal reparametrization, the inverse of the stereographic projection. This result can be seen as a variant of the Hopf's theorem where the hypothesis of conformality is replaced by a bound on the area.   
\begin{lemma}[lemma A.1 of \cite{BC2}]
\label{bcl}
Let $\omega \in L^{1}_{loc}(\R^2,\R^3)$ which satisfies 
\beq
\label{eqlim}
\begin{split}
&\Delta \omega = -2\, \omega_x \wedge \omega_y ,\\
& \int_{\R^2} \vert\nabla  \omega \vert^2 dz <+\infty .
\end{split}
\eeq
Then $\omega$ has precisely the form 
$$\omega(z) =\pi^{-1}_{N} \left(\frac{P(z)}{Q(z)} \right) + C, $$
where $N\in S^2$, $P$ and $Q$ are polynomial, C is a constant and $\pi_{N}$ is the stereographic projection from the north pole $N$. In addition
$$\int_{\R^2} \vert \nabla \omega \vert^2 = 8\pi k \hbox{ with }
 k=max\{ deg P, deg\ Q\},$$
 provided that $\frac{P}{Q}$ is irreducible.
 \end{lemma}
 It could be useful to remark that the gradient of such an $\omega$ satisfies the following formula
 $$\vert \nabla \omega \vert = \frac{2\sqrt{2} \vert P' Q - Q'P\vert}{\vert P\vert^2 +\vert Q\vert^2} . $$
 Then we define a special class of solutions which will be important in what follows: the spheres which are parametrized only once.
 \begin{defi} A solution $\omega$ of (\ref{eqlim}) is said to be simple if
$$\omega(z) =\pi^{-1}_{N} \left(\frac{P(z)}{Q(z)} \right) + C, $$
with  $\frac{P}{Q}$ is irreducible and $max\{ deg P, deg\ Q\}= 1$.
\end{defi}
In particular, if $\omega$ is a simple solution of (\ref{eqlim}), then we have 
\beq
\label{oi}
\left\vert \nabla   \omega^\eps (x) \right\vert = O \left(\frac{\lambda^\eps}{\vert x- a^\eps\vert^2 +(\lambda^\eps)^2} \right), 
\eeq
where  $\omega^\eps = \omega\left( \frac{\, .\, - a^\eps}{\lambda^\eps} \right)$, $a^\eps$ and $\lambda^\eps$ are respectively a sequence of points in $\R^2$ and  a sequence of positive numbers.\\

Finally we give a generalisation of the result of Brezis and Coron for solutions defined on the disk or the half-plane with appropriate boundary conditions.
\begin{lemma}
\label{leqlim}
Let $\Omega= \h$ or $\D$ and $u: \Omega \rightarrow \R^2\times \R_{+}$ be such that 
\be
\begin{split}
&\Delta u=-2 u_x \wedge u_y,\\
&\langle u_x,u_y\rangle= \vert u_x\vert -\vert u_y \vert =0,\\
& \Vert \nabla u \Vert_{2}< + \infty,\\
&u_{\vert \partial \Omega} \subset \R^2\times \{ 0\} 
\end{split}
\ee
and such that the angle between $u(\Omega)$ and $\R^2\times \{ 0\} $ when it is defined is right. Then
$$u=C+ \pi^{-1}\left( \frac{P}{Q}\right)$$
where $\pi$ is the stereographic projection and $P$ and $Q$ are two polynoms of  $\C[z]$. Moreover, $u_{\vert \partial \Omega}$ describes a circle of radius one.
\end{lemma}

{\it Proof of lemma \ref{leqlim}:}\\
 
First of all we can assume that $\Omega=\D$. Indeed let $\phi : \overline{\D}\setminus\{1\} \rightarrow \overline{\h}$ defined by
 $$\phi(z) =-i\frac{z+1}{z-1}.$$
It is well known that this application is a conformal isomorphism. Hence if $u:\h \rightarrow \R^2\times \R_{+}$ satisfies the hypothesis of the lemma, it is the same for $\tilde{u}: \overline{\D} \setminus\{1\} \rightarrow \R^2\times \R_{+}$ defined by $\tilde{u}=u\circ \phi$. But, since $\Vert \nabla \tilde{u}\Vert_2 =\Vert \nabla u\Vert_2<+\infty$, thanks to the regularity theory $\tilde{u}$ can be extended smoothly  at $1$.\\ 
 
Then we can extend $u$ to the whole plane, setting
 $$u(z) =-\left(\begin{array}{r}u^1 \\u^2 \\-u^3\end{array}\right) \left(\frac{1}{\overline{z}}\right) \hbox{ pour tout } z\in \R^2\setminus \D .$$
This extension is  $C^1$ and moreover satisfies the hypothesis of lemma \ref{bcl}, since the energy is simply doubled by this extension. This proves the first part of the theorem. Finally we easily remark that $u_{\vert \partial\Omega}$ describe a circle of $1$, remarking that $u$ is equal, up to a sign change, to its Gauss map. But our hypothesis forces to the Gauss map to be contained in a great circle on the boundary, which achieves the proof of the theorem.\hfill$\square$

\section{Proof of theorem}
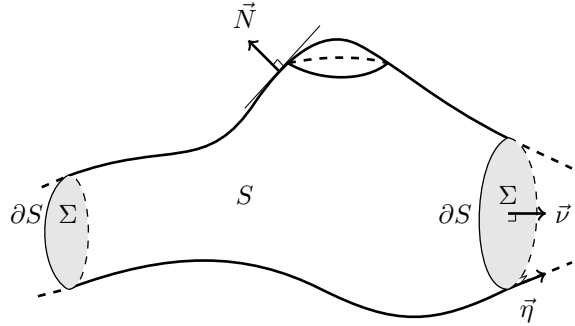
\begin{figure}[!h]\begin{center}\input{figbf}\end{center}\caption{ A surface with boundary.}\label{figbf}\end{figure}

The main idea is to apply the balancing formula to the boundary of our sequences of surfaces in order to detect the geometry of $\partial \Omega$. We remind that the balancing formula is an identity discovered by Kusner \cite{KKS} which concerns the shape of the boundary of a surface with constant mean curvature. Let $S$ be a surface with constant mean curvature. Then
\beq
\label{pbf3}
\int_{\partial S}  \vec{\eta}\, ds = 2 H_0 \int_{\Sigma}  \vec{\nu}\, d\sigma,
\eeq
where $\Sigma$ is a smooth surface having the same boundary than $S$, $\vec{\nu}$ is the normal of $\Sigma$ and  $\vec{\eta}$ is the conormal of $\partial S$, see figure \ref{figbf}. The reader will find a proof of (\ref{pbf3}) at chapter 7 of \cite{kenmotsu}.

In order to exploit this formula, we need a precise description of the behaviour of our surfaces. In order to obtain it, we start by proving that the sequence decomposes asymptoticly as a sum of spheres and hemispheres. Then we shall  prove that this decomposition contains at least one hemisphere, that is to say that the boundaries of our sequence of surfaces do not collapse to a point.\\  

Now we consider a smooth domain $\Omega$ of $\R^3$ and a sequence of embedded disks $\Sigma^\eps$ in $\Omega$ which satisfy the following assumptions
\begin{enumerate}[(i)]
\item  $\partial\Sigma^\eps \subset \partial\Omega$ and  $\partial\Sigma^\eps$ and $\partial\Omega$ meet orthogonally, 
\item  $\Sigma^\eps$ has \cmc\ equal to $\frac{1}{\eps}$,
\item the diameter and the area of $\Sigma^\eps$ are respectively a  $O(\eps)$ and a  $O(\eps^2)$.
\end{enumerate}
Up to translate $\Omega$ and to extract a subsequence of  $\Sigma^\eps$, we can assume that $\Sigma^\eps$ goes to $0$ and that $0 \in \partial \Sigma^\eps$. Then we rescale the space by a factor $\frac{1}{\eps}$ and we choose a conformal parametrization for our sequence of surfaces, that is to say a sequence of $u^\eps :\D \rightarrow \R^3$ such that 
\beq
\label{princ}
\begin{split}
&\begin{array}{cc}
\begin{cases}
\Delta u^\eps = -2 \, u_{x}^\eps \wedge u_{y}^\eps , \\ 
\langle u_x^\eps, u_y^\eps\rangle = \Vert u_x^\eps \Vert -\Vert u_y^\eps \Vert =0,
\end{cases}
& \hbox{ on } \D, 
\end{array}\\
& \Vert u^\eps\Vert_\infty  =O(1) \hbox{ and }  \Vert\nabla u^\eps\Vert_2  =O(1) ,\\
& u^\eps(\partial \D) \subset \partial \Omega^\eps \hbox{ and } \left\langle  u_{x}^\eps \wedge u_{y}^\eps, N^\eps \right\rangle =0 \hbox{ on } \partial \D,
\end{split}
\eeq
where $\Omega^\eps =\frac{1}{\eps}\Omega$ and $N^\eps$ is the exterior normal of $\partial\Omega^\eps$. The regularity of such a sequence of functions depends on the regularity of the surface where its free boundary lives. Here, since $\partial \Omega$ is smooth,  our sequence is smooth up to the boundary.
\subsection{Decomposition of $u^\eps$ as a sum of spheres and hemispheres.}
\label{sXX}
We start performing a decomposition of our surfaces as  a sum of spheres and hemispheres in the spirit of what has been done by Brezis and Coron in \cite{BC2}. However there are two big changes. On the one hand, there are two limit solutions (sphere and hemisphere). On the other hand, we must obtain an $L^\infty $-estimate rather than estimates on the gradient, while the equation lends itself much better to obtain estimates on the gradient.
\begin{thm}
\label{de'}
Let $u^\eps$ be a sequence of maps in $C^2(\overline{\D})$ which are  non-constant solutions of (\ref{princ}), then either $u^\eps$ converge uniformly to $0$ or there exist $p\in\N$ and
\begin{enumerate}[(i)]
\item   $\omega^1,\dots,\omega^p $  non-constant solutions of (\ref{eqlim}), 
\item  $ a^\eps_1, \dots, a^\eps_p$  sequences of $\overline{\D}$, and
\item  $\lambda^\eps_1 , \dots , \lambda^\eps_p$ sequences of positive real numbers such that $\ds \lim_{\eps \rightarrow 0} \lambda^\eps_i <+\infty $, 
\end{enumerate} 
such that, for a subsequence $u^\eps$ (still denote $u^\eps$), we get
\beq
\label{A'}
\tag{A}
 u^\eps_i \rightarrow \omega^i \hbox{ in } C^2_{loc}\left(\Omega_i\setminus S_i\right) \hbox{ as } \eps \rightarrow 0 \hbox{ for all } 1\leq i\leq p,
 \eeq
where $u^\eps_i= u^\eps (\lambda^\eps_i\ . + a_i^\eps )$, $\ds \Omega_i =\lim_{\eps \rightarrow 0} \left\{ z\in\R^2 \hbox{ s.t. } \lambda^\eps_i\ . + a_i^\eps\in \D \right\}$ and\\ $\ds S_i= \lim_{\eps \rightarrow 0}  \left\{ \frac{a_j^\eps-a_i^\eps}{\lambda_i^\eps}\hbox{ s.t. } j\in \{1,\dots,p\}\setminus \{ i\}\right\} $.
\beq
\label{ortho'}
\tag{B}
\lim_{\eps\rightarrow 0} \frac{d_{i}^\eps(a_j^\eps)}{\lambda^\eps_j} + \frac{d_{j}^\eps(a_i^\eps)}{\lambda^\eps_i} = +\infty  \hbox{ for all }i \not= j ,
\eeq
where $d_{i}^\eps(x)=\sqrt{(\lambda^\eps_i)^2 +\vert a_i^\eps -x\vert^2}$, 
\beq
\tag{C}
\label{idec}
d\left(u^\eps(\D), \cup_{i=1}^p B_i \right) \rightarrow 0  \hbox{ as } \eps\rightarrow 0,
\eeq 
where $B_i$ are the limit set of $\omega^\eps_i (\D)$ as $\eps$ goes to zero, with $\omega^\eps_i = \omega_i\left(\frac{\ .\ -a_i^\eps}{\lambda_i^\eps} \right)$,  that is to say some  spheres and hemispheres.
\end{thm}

{\it Proof of theorem \ref{de'} :}\\

We are going to extract the bubble by induction, the process will stop thanks to our uniform estimate on the energy of  $u^\eps$. In fact, such an extraction will be done until a "weak estimate" is not satisfied on the reminder, which is by now an almost classical technic since the work of Druet, Hebey and Robert about strong estimate for sequences of solution of Yamabe-type equation, see \cite{DruetHebeyRobert}. The advantage of this method is to insure a $C^2_{loc}$-converge rather than an $H^1$-converge.\\

{\bf Let $k\geq 1$, we say that $u^\eps$ satisfies the property $(P_k)$ if there exist }\\
\begin{enumerate}[(i)]
\item  $\omega^1,\dots,\omega^k $ non-constant solution of (\ref{eqlim}),
\item  $ a^\eps_1, \dots, a^\eps_k$  sequences of $\overline{\D}$ and
\item  $\lambda^\eps_1, \dots,\lambda^\eps_k$ sequences of  positive real numbers such that $\ds \lim_{\eps \rightarrow 0} \lambda^\eps_i < +\infty$ , 
\end{enumerate} 
such that, for a subsequence $u^\eps$ (still denoted $u^\eps$), we get
\beq
\tag{$A_k$}
\label{ui'}
 u^\eps_i \rightarrow \omega_i \hbox{ in } C^2_{loc}\left(\Omega_i\setminus S_i \right) \hbox{ as } \eps \rightarrow 0 \hbox{ for all } 1\leq i\leq k,
 \eeq
where $u^\eps_i= u^\eps (\lambda^\eps_i\ . + a_i^\eps )$, $\ds \Omega_i =\lim_{\eps \rightarrow 0} \left\{ z\in\R^2 \hbox{ s.t. } \lambda^\eps_i\ . + a_i^\eps\in \D \right\}$ and\\ $\ds S_i= \lim_{\eps \rightarrow 0}  \left\{ \frac{a_j^\eps-a_i^\eps}{\lambda_i^\eps}\hbox{ s.t. } j\in \{1,\dots,p\}\setminus \{ i\}\right\} $.
 \beq
 \tag{$B_k$}
\label{ortho1'}
 \frac{d_{i}^\eps(a_j^\eps)}{\lambda^\eps_j} + \frac{d_{j}^\eps(a_i^\eps)}{\lambda^\eps_i} \rightarrow +\infty  \hbox{ as } \eps \rightarrow 0 \hbox{ for all } i \not= j,
\eeq
where $d_{i}^\eps(x)=\sqrt{(\lambda^\eps_i)^2 +\vert a_i^\eps -x\vert^2}$. Moreover, when $\Omega_i \not= \R^2$, $\omega^i_{\vert\partial \Omega_i} $ describe (perhaps several times) a circle of radius  $1$.\\

{\bf Claim 1: If $(P_k)$ holds for some $k\geq 1$, then either $(P_{k+1})$ holds or 
\beq
\lim_{\eps \rightarrow 0} \sup_{z\in \D}  \left(\min_{1\leq i \leq k} d_i^\eps (z)\right) \left\vert \nabla \left(u^\eps -\sum_{i=1}^k  \omega^\eps_i\right) (z)\right\vert =0,
\eeq
where $\omega^\eps_i = \omega_i\left(\frac{\, . \, -a_i^\eps}{\lambda_i^\eps} \right)$ .
}\\

{\it Proof of claim $1$:}\\

Assume that $(P_k)$ holds and that there exists $\gamma_0>0$  and a subsequence $u^\eps$ (still denoted $u^\eps$) such that 
\beq
\label{gam0'}
 \sup_{z\in \D}  \left(\min_{1\leq i \leq k} d_i^\eps (z)\right) \left\vert \nabla \left(u^\eps-\sum_{i=1}^k  \omega^\eps_i\right) (z) \right\vert \geq \gamma_0.
\eeq
Let  $a^\eps_{k+1} \in \overline{\D}$ be such that 
\be
\left(\min_{1\leq i\leq k} d_i^\eps (a_{k+1}^\eps) \right) \left\vert \nabla  \left(u^\eps -\sum_{i=1}^k  \omega^\eps_i\right)  (a_{k+1}^\eps) \right\vert= \sup_{z\in \D} \left( \min_{1\leq i \leq k} d_i^\eps (z) \right) \left\vert \nabla  \left(u^\eps- \sum_{i=1}^k  \omega^\eps_i\right)  (z) \right\vert .
\ee 
We define $\lambda^\eps_{k+1}$ by 
\be
\left\vert \nabla\left(u^\eps-\sum_{i=1}^k  \omega^\eps_i\right) (a_{k+1}^\eps)\right\vert  = \frac{1}{\lambda_{k+1}^\eps} . 
\ee
Remarking that $\ds \min_{1\leq i\leq k} d_i^\eps (a_{k+1}^\eps)$ is bounded, it is clear that
\beq
\label{l'} 
\lim_{\eps\rightarrow 0}\lambda_{k+1}^\eps <+\infty .
\eeq
There are now two cases to consider.\\

{\bf First case :
\beq
\label{c1'}
 \lim_{\eps\rightarrow 0} \frac{\displaystyle \min_{1\leq i \leq k} d_i^\eps (a_{k+1}^\eps)}{\lambda_{k+1}^\eps}  = +\infty .
\eeq
}
In that case, ($B_{k+1}$) is automatically satisfied. We set
\be
u^\eps_{k+1} (z)= u^\eps \left( \lambda_{k+1}^\eps z  +a_{k+1}^\eps\right) \hbox{ for all } z\in \overline{\Omega}^\eps_{k+1}
\ee
where $\displaystyle \Omega^\eps_{k+1} =\left\{ z\in\R^2 \hbox{ s.t. }  \lambda_{k+1}^\eps z  +a_{k+1}^\eps\in \D\right\} $. Let $z\in \Omega^\eps_{k+1}$, we get
\beq
\label{c2'}
\begin{split}
\vert \nabla u^\eps_{k+1}(z)\vert &=  \lambda_{k+1}^\eps \vert \nabla u^\eps ( \lambda_{k+1}^\eps z + a_{k+1}^\eps)\vert \\
&\leq  \lambda_{k+1}^\eps \left\vert \nabla \left(u^\eps - \sum_{i=1}^k  \omega^\eps_i\right)(  \lambda_{k+1}^\eps z + a_{k+1}^\eps)\right\vert \\
&+  \lambda_{k+1}^\eps \left\vert \nabla \left(\sum_{i=1}^k  \omega^\eps_i\right)( \lambda_{k+1}^\eps z + a_{k+1}^\eps)\right\vert .
\end{split}
\eeq
Thanks to (\ref{ui'}) and (\ref{c1'}), we easily get that 
\beq
\label{c3'}
\lambda_{k+1}^\eps \left\vert \nabla \left(\sum_{i=1}^k  \omega^\eps_i\right)(\lambda_{k+1}^\eps z +a_{k+1}^\eps )\right\vert =o(1) .
\eeq
Then, using the definition of $a_{k+1}^\eps$, (\ref{c1'}), (\ref{c2'}) and (\ref{c3'}), we get
\beq
\label{c4'}
\vert \nabla u^\eps_{k+1}(z)\vert  \leq \frac{\min_{1\leq i\leq k} d_i^\eps (a_{k+1}^\eps)}{\min_{1\leq i\leq k} d_i^\eps (\lambda_{k+1}^\eps z+ a_{k+1}^\eps)} +o(1)=1+o(1) .
 \eeq
Then $\vert \nabla u^\eps_{k+1}\vert$ is bounded on every compact subset of $\overline{\Omega}^\eps_{k+1}$. Moreover, thanks to conformal invariance of our equation, $u^\eps_{k+1}$ still satisfies (\ref{princ}). Hence, using standard elliptic theory, see section \ref{sregul} and \cite{GT}, we see that there exist a subsequence of $u^\eps$ (still denoted  $u^\eps$) and $\omega^{k+1}\in C^2(\overline{\Omega}_{k+1})$ such that  
\be
u^\eps_{k+1} \rightarrow \omega^{k+1} \hbox{ in } C^2_{loc} \left(\overline{\Omega_{k+1}}\right) \ee
and 
\be
\Delta \omega^{k+1} = -2\,  \omega^{k+1}_x \wedge \omega^{k+1}_y \hbox{ on } \Omega_{k+1},
\ee
where $\displaystyle\Omega_{k+1} =\lim_{\eps \rightarrow 0} \Omega_{k+1}^\eps$. Here there are again two cases: either $\Omega_{k+1}$ is the whole plane or this is a disk or an half-plane (which is conformally equivalent).  In the case of the disk or the half-plane, the boundary condition pass to the limit, that is to say, up to rotation,
\be
\begin{split}
&\omega^{k+1}(\partial \Omega_{k+1}) \subset \R^2\times\{0\} \\
&\hbox { and }\\
& \left\langle  \omega^{k+1}_{x} \wedge \omega^{k+1}_{y}, N \right\rangle =0 \hbox{ on } \partial \Omega_{k+1} ,
 \end{split}
 \ee
where $N=(0,0,1)$. Moreover, thanks to conformal invariance of $\Vert\nabla .\, \Vert_2$, up to extract a subsequence, we get
$$u^\eps_{k+1}\rightharpoonup \omega^{k+1} \hbox{ in } L^2 (\R^2) $$
and
$$\Vert \nabla \omega^{k+1} \Vert_2 \leq \liminf_{\eps\rightarrow 0} \Vert \nabla u^\eps_{k+1}\Vert_2= \liminf_{\eps\rightarrow 0}\Vert \nabla u^\eps\Vert_2<+\infty .$$
Finally, thanks to lemmas \ref{bcl} and  \ref{leqlim}, $\omega^{k+1}$ has the desired shape, moreover $\omega^{k+1}$ is non-constant since $\vert \nabla \omega^{k+1}(0)\vert=1$. This achieves the proof of  $(P_{k+1})$ in the first case.\\

{\bf Second case :
\beq
\label{29'}
 \lim_{\eps\rightarrow 0} \frac{\ds \min_{1\leq i \leq k} d_i^\eps (a_{k+1}^\eps)}{\lambda_{k+1}^\eps}  = \gamma >0 .
\eeq
}

First of all, we need to show that ($B_{k+1}$) holds. We assume by contradiction that ($B_{k+1}$) does not hold, then, up to extract a subsequence, there exists $1\leq i_0 \leq k$ such that  
\beq
\label{c12'} 
d_{k+1}^\eps (a_{i_0}^\eps)=O(\lambda_{i_0}^\eps) \hbox{ and } d_{i_0}^\eps (a_{k+1}^\eps) =O(  \lambda_{k+1}^\eps) . 
\eeq
On the one hand, (\ref{c12'}) gives
\beq
\label{c13'}
\begin{split}
\lim_{\eps\rightarrow 0 }\frac{\lambda_{k+1}^\eps}{\lambda_{i_0}^\eps}= c 
\hbox{ and }
\vert a_{i_0}^\eps - a_{k+1}^\eps \vert = O(\lambda_{i_0}^\eps),
\end{split}
\eeq
where $c$ is a positive constant. On the other hand, thanks to (\ref{ui'}) and (\ref{ortho1'}), we get 
\beq
\label{c14'}
 \nabla \left(\left(u^\eps-\sum_{i=1}^k  \omega^\eps_i\right) (\lambda_{i_0}^\eps\ . + a_{i_0}^\eps)\right) \rightarrow 0 \hbox{ in }C^2_{loc} \left(\Omega_{i_0} \setminus S_{i_0} \right). 
 \eeq
Hence, thanks to (\ref{29'}) and (\ref{c13'}), we necessarily get    
$$ d\left(\frac{a_{k+1}^\eps - a_{i_0}^\eps}{\lambda_{i_0}^\eps}, S_{i_0}\right)= o(1).$$ 
Let $j\in \{1,\dots,k\}\setminus\{i_0\}$ be such that 
$$\left\vert \frac{a_{k+1}^\eps - a_{j}^\eps}{\lambda_{i_0}^\eps}\right\vert= o(1).$$ 
Using (\ref{29'}) et (\ref{c13'}), we remark that for  $\eps$ small enough, 
$$ \frac{\lambda_j^\eps}{\lambda_{k+1}^\eps} \geq \frac{\gamma}{2} ,$$
and, using again (\ref{c13'}), we remark that, for $\eps$ small enough,  
$$
\frac{\lambda_j^\eps}{\lambda_{i_0}^\eps}  \geq \frac{\gamma}{4c} ,
$$
Since $\frac{a_{i_0}^\eps -a_j^\eps}{\lambda_{i_0}^\eps} =O(1)$ and $i_0$ and $j$ satisfy (\ref{ortho1'}), we necessarily get 
$$
\lambda_{i_0}^\eps =o(\lambda_j^\eps) .
$$
Hence, for all $j$ such that $\frac{a_{k+1}^\eps -a_j^\eps}{\lambda_{i_0}^\eps} =o(1)$, we get
$$
\lambda_{i_0}^\eps =o(\lambda_j^\eps) .
$$
In particular, thanks to (\ref{ui'}), there exits $\delta >0$ such that for all $z\in B(0,\delta)$, we have 
$$\lambda_{i_0}^\eps \vert \nabla \omega_i^\eps (a_{k+1}^\eps + z \lambda_{i_0}^\eps)\vert =o(1) \hbox{ for all } i\not= i_0 .$$
We easily see that
$$\lambda_{i_0}^\eps \vert \nabla u^\eps \vert = O(1) \hbox{ on } B(a_{k+1}^\eps, \delta \lambda_{i_0}^\eps) .$$
Using standard elliptic theory, up to extract a subsequence, we get that $u^\eps_{i_0}$ converge to $\omega^{i_0}$ in $\ds C^2_{loc} \left(B\left(c_{k+1},\frac{\delta}{2}\right)\right)$ where $\ds c_{k+1}=\lim_{\eps\rightarrow 0} \frac{a_{k+1}^\eps - a_{i_0}^\eps}{\lambda_{i_0}^\eps}$. Then we deduce that $$ \vert \nabla (u^\eps_{i_0} -\omega^{i_0})(a^\eps_{k+1})\vert \rightarrow 0,$$
which leads to
$$
\lambda_{i_0}^\eps \left\vert\nabla \left(\left(u^\eps-\sum_{i=1}^k  \omega^\eps_i\right) (a_{k+1}^\eps)\right)\right\vert \rightarrow 0, 
 $$
which, thanks to (\ref{c14'}), is  a contradiction with (\ref{29'}) and proves ($B_{k+1}$).\\ 

Now, we set
$$u^\eps_{k+1} = u^\eps(\lambda_{k+1}^\eps\ .+ a_{k+1}^\eps) \hbox{ for all } z\in\Omega^\eps_{k+1} ,$$
where $\displaystyle \Omega^\eps_{k+1} =\left\{ z\in\R^2 \hbox{ s.t. }  \lambda_{k+1}^\eps z  +a_{k+1}^\eps\in \D\right\} $. Let $z\in  \overline{\Omega^\eps_{k+1}} \setminus \{S_{k+1}\}$, we have
\beq
\label{c22'}
\begin{split}
\vert \nabla u^\eps_{k+1}(z)\vert &=  \lambda_{k+1}^\eps \vert \nabla u^\eps (\lambda_{k+1}^\eps z+ a_{k+1}^\eps)\vert\\
&\leq  \lambda_{k+1}^\eps \left\vert \nabla \left(u^\eps - \sum_{i=1}^k  \omega^\eps_i\right)(\lambda_{k+1}^\eps z+ a_{k+1}^\eps)\right\vert \\
&+  \lambda_{k+1}^\eps \left\vert \nabla \left(\sum_{i=1}^k  \omega^\eps_i\right)(\lambda_{k+1}^\eps z + a_{k+1}^\eps)\right\vert .
\end{split}
\eeq
Thanks to (\ref{ui'}) and (\ref{29'}), we obtain
\beq
\label{c32'}
\lambda_{k+1}^\eps \left\vert \nabla \left(\sum_{i=1}^k  \omega^\eps_i\right)(\lambda_{k+1}^\eps\ .+ a_{k+1}^\eps)\right\vert =O\left(\frac{1}{d(z,S_{k+1})} + \vert z \vert +1 \right).
\eeq
With the concention that $d(z,\emptyset)=+\infty$.\\

Then using the definition of $a_{k+1}^\eps$, (\ref{c22'}) and (\ref{c32'}), we get
\beq
\label{c33'}
\begin{split}
\vert \nabla u^\eps_{k+1}(z)\vert  &\leq \frac{\min_{1\leq i\leq k} d_i^\eps (a_{k+1}^\eps)}{\min_{1\leq i\leq k} d_i^\eps (\lambda_{k+1}^\eps z + a_{k+1}^\eps)} +  O\left(\frac{1}{d(z,S_{k+1})} + \vert z \vert +1 \right)\\
&=O\left(\frac{1}{d(z,S_{k+1})} + \vert z \vert +1 \right) .
\end{split}
 \eeq
Then $\vert \nabla u^\eps_{k+1}\vert$ is bounded on every compact subset of $\overline{\Omega^\eps_{k+1}} \setminus S_{k+1}$. Moreover, thanks to the conformal invariance of our equation, $u^\eps_{k+1}$ still satisfies (\ref{princ}). Hence, thanks to the standard elliptic theory, see section \ref{sregul} and \cite{GT}, there exists a subsequence of $u^\eps$ (still denoted $u^\eps)$ and $\omega^{k+1}\in C^2(\overline{\Omega_{k+1}} \setminus S_{k+1})$  such that 
$$u^\eps_{k+1} \rightarrow \omega^{k+1} \hbox{ in } C^1_{loc} (\overline{ \Omega_{k+1}}\setminus S_{k+1}) ,$$
and
\be
\Delta \omega^{k+1}  = -2\,  \omega^{k+1}_x \wedge \omega^{k+1}_y \hbox{ on } \Omega_{k+1}\setminus S_{k+1}
\ee
where $\displaystyle\Omega_{k+1} =\lim_{\eps \rightarrow 0} \Omega^\eps_{k+1}$. As before, there are two possibilities; either  $\Omega_{k+1}$ is the whole plane or this a disk or an half-plane (which is conformally equivalent). If it is a disk or an half-plane, the boundary condition pass to the limit, that is to say, up to rotation,
\be
\begin{split}
&\omega^{k+1}(\partial \Omega_{k+1}) \subset \R^2\times\{0\} \\
&\hbox { et }\\
& \left\langle  \omega^{k+1}_{x} \wedge \omega^{k+1}_{y}, N \right\rangle =0 \hbox{ on } \partial \Omega_{k+1} ,
 \end{split}
 \ee
where $N=(0,0,1)$.\\
Moreover, thanks the conformal invariance of $\Vert\nabla .\, \Vert_2$, up to extract a subsequence, we get 
$$u^\eps_{k+1}\rightharpoonup \omega^{k+1} \hbox{ in } L^2 (\R^2) $$
and
$$\Vert \nabla \omega^{k+1} \Vert_2 \leq \liminf_{\eps\rightarrow 0} \Vert \nabla u^\eps_{k+1}\Vert_2= \liminf_{\eps\rightarrow 0}\Vert \nabla u^\eps\Vert_2<+\infty .$$
Then, $\omega^{k+1}$ is a solution of (\ref{eqlim}) on $\Omega_{k+1}$, and $\omega^{k+1}$ has the desired shape. Finally, we need to show that $\omega^{k+1}$ is non-constant. This is trivial if  $0\not\in S_{k+1}$, since in that case $\vert \nabla \omega^{k+1}(0)\vert=1$. Else, for all $i_0$ such that 
$$\frac{\vert a_{i_0}^\eps - a_{k+1}^\eps\vert}{\lambda_{k+1}^\eps} =o(1),$$
thanks to (\ref{29'}) and ($B_{k+1}$), we get
$$ \lambda_{i_0}^\eps =o(\lambda_{k+1}^\eps).$$
Then mimicking the argument of the proof of ($B_{k+1}$) we show that
$$ \nabla u^\eps_{k+1} \rightarrow \nabla \omega^{k+1} \hbox{ on } B(0, \delta),$$
where $\delta>0$. This leads in every case to $\vert \nabla \omega^{k+1}(0)\vert=1$ and then  $\omega^{k+1}$ is non-constant. This proves ($P_{k+1}$) in this second case. The study of these two cases ends the proof of claim 1.\hfill$\square$\\

Then, we need to prove a claim about the energy of a sum of bubbles. In fact, using ($B_k$), we show that the bubbles do not interact in a weak sense and that each one provides at least  the energy of a simple hemisphere,that is to say $4\pi$. \\ 

{\bf Claim 2: Let $k\in \N^*$ and
\begin{enumerate}[(i)]
\item   $\omega^1,\dots,\omega^k $  non-constant solutions of (\ref{eqlim}), 
\item   $ a^\eps_1, \dots, a^\eps_k $ sequences of $\overline{\D}$, and
\item  $\lambda^\eps_1, \dots,\lambda^\eps_k,$ sequences of positive real numbers such that  $\ds \lim_{\eps \rightarrow 0} \lambda^\eps_i =0$, 
\end{enumerate} 
such that, with $u^\eps$, they satisfy ($P_k$). Then
$$ \liminf_{\eps \rightarrow 0} \Vert \nabla u^\eps \Vert_2^2 \geq \sum_{i=1}^k \Vert \nabla \omega^i \Vert^2_2 \geq 4\pi k .$$
}\\

{\it Proof of claim 2 :}\\

Let $R$ be a positive real number, thanks to (\ref{ortho1'}), for $\eps$ small enough, we get
$$
\int_{\D} \vert \nabla u^\eps \vert^2   dz \geq \sum_{i=1}^k
\int_{\D\cap B(a^\eps_i, R\lambda_i^\eps)\setminus \Omega_i^\eps (R)} \vert \nabla u^\eps\vert^2 dz,
$$
where 

$$\Omega_i^\eps (R) = \left\{z\in B(a^\eps_j, R \lambda_j  ^\eps)  \hbox{ where } j \hbox{ is such that} \lim_{\eps \rightarrow 0} \frac{\lambda_j^\eps}{\lambda_i^\eps}=0\right\}.$$ 
Then, thanks to ($A_k$), we get
\beq
\begin{split}
\int_{\R^2} \vert \nabla u^\eps \vert^2   dz
 &\geq   \sum_{i=1}^k \int_{\left(B\left(-\frac{a_i^\eps}{\lambda_i^\eps},\frac{1}{\lambda_i^\eps}\right) \cap B(0, R)\right)\setminus \Omega_i (R)} \vert \nabla \omega_i \vert^2 dz + \delta_{\eps,R} \\
 & \geq 4\pi k + \delta_{\eps,R}
\end{split}
\eeq
where $\ds \Omega_i (R)= \cup_{x\in S_i} B(x,\frac{1}{R})$ and $\ds \lim_{R\rightarrow +\infty } \lim_{\eps \rightarrow 0} \delta_{\eps,R}=0$. Here we bound the energy of a solution by the smallest possible, that is to say the energy of an hemisphere. \hfill$\square$\\

{\bf Proof of the theorem:}\\

We start setting $a_1^\eps  \in \overline{\D}$ and $\lambda_1^\eps$ as
\be \left\vert \nabla  u^\eps  (a_{1}^\eps) \right\vert= \sup_{z\in \D}  \left\vert \nabla u^\eps  (z) \right\vert 
\ee 
and
$$\left\vert \nabla u^\eps  (a_{1}^\eps)\right\vert  = \frac{1}{\lambda_{1}^\eps} . $$
Either  $\lambda_1^\eps$ goes to infinity and then  $u^\eps$ converges uniformly to $0$ which prove the theorem. Or we set
$$u_1^\eps (z)= u^\eps(a_1^\eps +\lambda_1^\eps z )\hbox{ for all } z\in \Omega_1^\eps$$
where $\ds \Omega_1^\eps=\left\{ z\in R^2 \hbox{ s.t. } a_1^\eps +\lambda_1^\eps z\in \D\right\}$.\\

It is clear that   $\vert \nabla u^\eps_{1}\vert$ is bounded on every compact subset of $\overline{\Omega^\eps_1}$. Moreover thanks to conformal invariance of our equation, $u^\eps_{1}$ still satisfies (\ref{princ}). Hence, applying standard elliptic theory, see section \ref{sregul} et \cite{GT}, We see that there exists a subsequence of $u_1^\eps$ (still denoted  $u_1^\eps)$ and $\omega^{1}\in C^2(\Omega_1)$ such that  
$$u^\eps_{1} \rightarrow \omega^{1} \hbox{ in } C^2_{loc} (\overline{\Omega_{1}}) $$
et 
\be
\Delta \omega^{1} = -2\,  \omega^{1}_x \wedge \omega^{1}_y \hbox{ on } \Omega_1 ,
\ee
where $\displaystyle\Omega_1 =\lim_{\eps \rightarrow 0} \Omega^\eps_1$. There are two possibilities for $\Omega_1$; it is either the whole plane or disk  or an half-plane (which is conformally equivalent). In the last case, the boundary condition passes to the limit, that is to say, up to rotation,
\be
\begin{split}
&\omega^{1}(\partial \Omega_{1}) \subset \R^2\times\{0\} \\
&\hbox { and }\\
& \left\langle  \omega^{1}_{x} \wedge \omega^{1}_{y}, N \right\rangle =0 \hbox{ on } \partial \Omega_{k} ,
 \end{split}
 \ee
where $N=(0,0,1)$. Moreover, thanks to conformal invariance of $\Vert\nabla .\, \Vert_2$, up to extract a subsequence, we get
$$u^\eps_{1}\rightharpoonup \omega^{1} \hbox{ in } L^2 (\R^2) $$
et
$$\Vert \nabla \omega^{1} \Vert_2 \leq \liminf_{\eps\rightarrow 0} \Vert \nabla u^\eps_{1}\Vert_2= \liminf_{\eps\rightarrow 0}\Vert \nabla u^\eps\Vert_2<+\infty .$$
Then, thanks to lemmas \ref{bcl} and \ref{leqlim},  $\omega^{1}$ has the desired shape. Finally $\omega^{1}$ is non-constant since $\vert \nabla \omega^{1}(0)\vert=1$.\\

Now we can start our induction. Indeed, thanks to claim 1 and 2 and the fact that the energy is uniformly bounded, there exists $k\in \N^*$ such that  ($P_k$) is satisfies and 
\beq
\label{ii'}
\lim_{\eps \rightarrow 0} \sup_{z\in \D}  \left(\min_{0\leq i\leq k} d_i^\eps (z)\right) \left\vert \nabla \left(u^\eps -\sum_{i=0}^k  \omega^\eps_i\right) (z)\right\vert =0,
\eeq
where $\ds\omega^\eps_i = \omega_i\left(\frac{\, . \, - a_i^\eps}{\lambda_i^\eps} \right)$. This proves (\ref{A'}) and (\ref{ortho'}).\\

It suffices to show (\ref{idec}) to conclude. We start with the following claim.\\

 {\bf Claim 3 : 
 \beq
\label{h1dec'}
\left\Vert\nabla \left( u^\eps -\sum_{i=1}^k \omega_i^\eps \right)\right\Vert_2 \rightarrow 0 \hbox{ when } \eps \rightarrow 0.
\eeq
}\\
 
 {\it Proof of claim $3$ :}\\
 
We set
$$R^\eps =u^\eps -\sum_{i=1}^k \omega_i^\eps$$ 
and we assume that there exists $\delta>0$ such that 
$$\Vert \nabla R^\eps \Vert_2 \geq \delta .$$
With those assumptions, we are going to prove the existence of a new bubble which will contradict (\ref{ii'}). In order to find this bubble, we follow the method developped in  \cite{BC2}.\\

First of all, we introduce the concentration function
$$C^\eps (t) = \sup_{z\in \D} \int_{B(z,t)} \vert \nabla R^\eps \vert^2\, dz.$$
It is clear that $C^\eps$ is continous, increasing with respect to $t$ and that $C^\eps(0)=0$. We fix $\nu$ such that 
$$0 < \nu < \min \{\frac{1}{2C_0}, \frac{\delta}{2} \} ,$$
where $C_0$ is the constant in the Wente inequality given by lemma \ref{Xlem}. Hence there exists $a^\eps \in \overline{\D}$ 
and  $\lambda^\eps >0$ such that   
 $$C^\eps (\lambda^\eps) =  \int_{B(a^\eps,\lambda^\eps)} \vert \nabla R^\eps \vert^2\, dz= \nu.$$
Then we rescale around $a^\eps$, setting $\tilde{f} =f(\lambda^\eps\, . +a^\eps)$ for all $\ds z\in \Omega^\eps =\left\{ z\in\R^2 \hbox{ s.t. } \lambda^\eps z + a^\eps \in \D \right\}$, and we get
$$\int_{\Omega^\eps} \vert \nabla \tilde{R}^\eps \vert^2 dz = \Vert \nabla R^\eps \Vert_2^2 \leq C, $$
and
$$\Vert \tilde{R}^\eps \Vert_\infty \leq C, $$
where $C$ is a positive real. Moreover, using (\ref{princ}), we remark that $\tilde{R}^\eps$ satisfy
\be
\Delta \tilde{R}^\eps = -2\, \tilde{R}^\eps_x \wedge \tilde{R}^\eps_y +O \left(\sum_{i=0}^k \vert \nabla \tilde{\omega}_i^\eps\vert \left( \sum_{ j\not= i}  \vert \nabla \tilde{\omega}_j^\eps\vert  + \vert \nabla \tilde{R}^\eps\vert\right)  \right)
\ee
However, thanks to (\ref{ortho'}), we get
$$\vert \nabla \tilde{\omega}_i^\eps\vert \vert \nabla \tilde{\omega}_j^\eps\vert \rightarrow 0 \hbox{ in } L^1_{loc}(\overline{\Omega^0}) \hbox{ pour } i\not= j$$
and, thanks to (\ref{ii'}), we get
 $$\vert \nabla \tilde{\omega}_i^\eps\vert \vert \nabla \tilde{R}^\eps\vert  \rightarrow 0 \hbox{ in } L^1_{loc}(\overline{\Omega^0}) \hbox{ for all } i,$$
with $\ds \Omega^0 = \lim_{\eps\rightarrow 0} \Omega^\eps$. Finally
$$ \Delta \tilde{R}^\eps= -2\, \tilde{R}^\eps_x \wedge \tilde{R}^\eps_y +h^\eps ,$$
where $h^\eps\rightarrow 0$ in $L^1_{loc}(\overline{ \Omega^0}) $ when $\eps \rightarrow 0$. Then, up to extract a subsequence, we get
$$ \tilde{R}^\eps \rightarrow R  \hbox{ p.p. on } \Omega^0$$
and
$$\nabla \tilde{R}^\eps \rightharpoonup \nabla R \hbox{ weakly in } L^2(\Omega^0). $$
Moreover $R$ is a weak solution of 
 $$\Delta R = -2 R_x \wedge R_y \hbox{ on } \Omega^0.$$ 
Now, thanks to our choice of $\nu$, we are going to show that the weak convergence is in fact strong. Let $v^\eps = \tilde{R}^\eps -R$, then $v^\eps$ satisfy 
\be \Delta v^\eps= -2\, v^\eps_x \wedge v^\eps_y -2 (v^\eps_x \wedge R_y + R_x \wedge v^\eps_y) +h^\eps .\ee
Moreover, thanks to corollary \ref{wente2}, there exists $\psi_\eps\in H^1_0(\Omega^0)$  a solution of
\be
\Delta \psi^\eps = -2 (v^\eps_x\wedge R_y+ R_x \wedge v^\eps_y)
\ee
satisfying
\beq
\label{h1'}
\Vert \nabla \psi^\eps  \Vert_2 + \Vert \psi^\eps  \Vert_\infty  \leq \Vert \nabla v^\eps  \Vert_2 \Vert \nabla R  \Vert_2  .
\eeq
However,
\be
\int_{\Omega^0} \vert \nabla \psi^\eps \vert^2 dz = -2 \int_{\Omega^0}  \langle \psi^\eps , v^\eps_x\wedge R_y+ R_x \wedge v^\eps_y\rangle\, dz .
\ee
Then using (\ref{h1'}),  we get that $\psi^\eps \wedge R_x$ and  $\psi^\eps \wedge R_y$  are bounded in  $L^2(\Omega^0)$. Hence, since  $\nabla v^\eps\rightarrow 0$ weakly in $L^2(\Omega^0)$, we get that
$$\int_{\Omega^0} \vert \nabla \psi^\eps \vert^2 dz \rightarrow 0. $$
Then we deduce that 
$$ \Delta v^\eps = -2\, v^\eps_x \wedge v^\eps_y  + g^\eps , $$
where $g^\eps\rightarrow 0$ in  $D'(\overline{\Omega^0})$.\\ 

Finally, let $\phi \in C^\infty_c(\overline{\Omega^0})$ be such that  $\operatorname{supp}(\phi)$ is  contained in a ball of radius $1$, using  lemma \ref{Xlem}, we get 
\be
\begin{split}
\int_{\Omega^0} \vert \nabla (\phi v^\eps) \vert^2 dz &= - 2 \int_{\Omega^0}  \langle  v^\eps , \phi  v^\eps_x\wedge  \phi v^\eps_y \rangle dz  + o(1) ,\\
& \leq 2\left(C_0\Vert \nabla  v^\eps_{| supp(\phi)} \Vert_2 \right) \Vert\nabla ( \phi  v^\eps)\Vert^2_2 +o(1). 
\end{split}
\ee
Thanks to our choice of   $\lambda^\eps$, we get $C_0 \Vert \nabla v^\eps_{| supp(\phi)} \Vert_2 \leq \frac{1}{2}$, which gives finally
$$ \int_{\Omega^0} \vert \nabla (\phi v^\eps) \vert^2 dz = o(1)$$
which proves
$$\nabla \tilde{R}^\eps \rightarrow \nabla R \hbox{ strongly in } L^2_{loc}(\overline{\Omega^0}) .$$

Indeed, we can remark that  $R$ isn't constant since $\Vert \nabla R \Vert_2 =\nu>0$. But, thanks to (\ref{ii'}), we have, for all  $z\in \R^2$, that  there exists $i$ such that 
\be
\vert \nabla \tilde{R}^\eps (z)\vert =o\left( \frac{1}{\sqrt{ \left(\frac{\lambda_{i}^\eps}{\lambda^\eps}\right)^2 + \left\vert z +\frac{a^\eps -a^\eps_i}{\lambda^\eps}\right\vert^2} }\right),
\ee
which is a contradiction and proves (\ref{h1dec'}).\hfill$\square$\\
 
In order to conclude, we have to transform this $H^1$-estimate in a $L^\infty$-estimate. An idea could be to use the Wente inequality as it is done by Brezis and Coron in \cite{BC2} in order to get $L^\infty$-estimate. But contrary to Brezis and Coron, here we don't control what happens on the boundary. In order to overpass this difficulty we are going to extend our surfaces.\\
Usually extend a surface across its boundary smoothly is not an easy fact, but here, thanks to the fact that our surfaces and the boundary of our domains meet orthogonally, this will be possible without perturbing too much the condition to be with constant mean curvature. 
\\
The idea is to reflect our surface through $\partial\Omega^\eps$ which is almost a plane so that our transformation will be almost an isometry (in fact a symmetry) and will almost conserve the mean curvature. Moreover the new surfaces will be at least $C^{1,1}$ thanks to the fact that our surfaces meet $\partial\Omega^\eps$ orthogonally.\\

Since $\partial\Omega^\eps$  converges uniformly to a plane , there exist a diffeomorphism $\psi^\eps :B(0,2R) \rightarrow \R^3$, where $R$ is chosen such that $u^\eps(\D)\subset B(0,R)$, which sends $\partial\Omega^\eps\cap B(0,2)$ in $\R^2\times \{0\}$ and which preserves the orthogonality on $\partial\Omega^\eps$. In fact it suffices to straighten up the local foliation of the normal bundle of $\partial\Omega^\eps$ to $\R^2\times \R$. Then now we get new surfaces which have almost constant mean curvature equal to $1$. Then we  extend our map to $S^2$. Here $S^2$ will be identified with the Riemann sphere $\hat{\C}$. We set
\be v^\eps (z)= s\left(v^\eps\left(\frac{1}{\overline{z}}\right)\right) \hbox{ for all } z\in \hat{\C}\setminus\D,\ee
where $s$ is the symmetry through $\R^2\times \{0\}$ and $v^\eps =\psi^\eps\circ u^\eps$. Using the fact that $v^\eps(\overline{\D})$ and $\R^2\times \{0\}$ meet orthogonally we easily show that $v^\eps$ is $C^{1,1}$.Then we set $\tilde{u}^\eps= \psi^{-1}\circ v^\eps$ which is also $C^{1,1}$ and its mean curvature uniformly converges to $1$.\\

Then we are in position to prove our theorem. Assume by contradiction that 
$$d\left(\tilde{\Sigma}^\eps, \cup_{i=1}^{k} \tilde{B}_i\right)\not\rightarrow 0,$$
where $\tilde{\Sigma}^\eps=\tilde{u}^\eps(S^2)$ and $\tilde{B}_i$ is the union of $B_i$ and its symmetry through $T_0\partial\Omega^\eps$.\\

Then there exits $y^\eps\in \tilde{\Sigma}^\eps$ such that
$$d\left(y^\eps, \cup_{i=1}^{k} \tilde{B}_i\right)\not\rightarrow 0.$$
Let $z^\eps \in \hat{C}$ be such that $\tilde{u}^\eps(z^\eps)=y^\eps$. We are going to prove that there is some area in a neighbourhood of $z^\eps$.  The idea is that if  a surface has bounded mean curvature, passes through the center of a ball, and has no boundary inside the ball, then it has to use a certain amount of area to leave the ball. Since the mean curvature is bounded, then  the Gaussian curvature of our surface is uniformly  bounded from above by a constant $K_0$. Let $r_0>0$ be such that 
$$ B\left(y^\eps,r_0\right) \cap \left(\cup_{i=1}^{k} \tilde{B}_i\right) =\emptyset.$$
Then using a Bishop comparison, like theorem III.4.2 of \cite{Chavel06}, we see that 

$$\hbox{Vol}\left( B\left(y^\eps,\frac{r_0}{2}\right) \cap \tilde{\Sigma}^\eps\right)\geq \hbox{Vol}\left(B_{\tilde{\Sigma}^\eps}\left(y^\eps,\frac{r_0}{2}\right)\right) \geq  \hbox{Vol}\left(B_{M_{K_0}}\left(y^\eps,\frac{r_0}{2}\right)\right) \geq C_0 r_0^2 ,$$
 where $M_{K_0}$ is the space of constant curvature $K^0$ and $C^0$ a positive constant. Hence we see that $u^\eps$  necessarily get some area in a neighbourhood of $z^\eps$ whose image is far from the bubbles, which is a contradiction with claim 3 since all the area of $u^\eps$ is devoted to cover the bubbles. This proves (C) and achieves the proof of the theorem.\hfill$\blacksquare$ \\

\subsection{There is at least one hemisphere in the decomposition}
In order to show our result, we have to eliminate  the sequence of surfaces whose boundaries collapse. It suffices to show that in theorem \ref{de'}, there is at least one bubble whose domain of definition is not the whole plane, that is to say there is at least one hemisphere.\\

Since our surfaces are embedded, we can assume that our bubbles are simple. In fact, we just need to prove that $max\{ deg P_i, deg\ Q_i\}= 1$ for all $1 \leq i\leq p$, with  $\omega_i =\pi_{P_i}\left( \frac{P_i}{Q_i}\right) $ where $\frac{P_i}{Q_i}$ is irreducible. But this is an easy consequence of the fact that our surfaces are embedded, (\ref{A'}) and the following lemma.
\begin{lemma}
\label{simple}
Let $u^\eps :B(0,1) \rightarrow \R^3$ a sequence of smooth embedding such that there exists $u^0\in C^1(B(0,1),\R^3 )$ and
\be
u^\eps \rightarrow u^0 \hbox{ in } C^{2}_{loc} (B(0,1)\setminus \{ 0\}) .
\ee
Then $u^0 $ can't be a multiple parametrization, that is to say there is no embedded $U_0\in C^1(B(0,1),\R^3 )$,  $\Phi \in \mathcal{O}(B(0,1),\C )$ an holomorphic function and an integer $k\geq 2$ such that
\be
\begin{split}
&u^0 = U^0 \circ \Phi\\
& \hbox{ and } \\
 &\Phi(z) = z^k+ o(\vert z\vert^k) \hbox{ as } z\rightarrow 0.
 \end{split}
 \ee
\end{lemma}

{\it Proof of the lemma \ref{simple} :} \\

First of all, up to a diffeomorphism of a neighbourhood of $0$, we can assume that 
\be
u^\eps \rightarrow U^0(z^l) \hbox{ in } C^{2}_{loc}(B(0,\delta)\setminus\{0\}) .
\ee
where $l\geq 2$ and $\delta>0$. Let $A_\delta = B \left(0,\frac{\delta}{2}\right)\setminus B\left(0,\frac{\delta}{3}\right)$ and $C_r$ be the cylinder of center $U^0(0)$, radius $r$ and orthogonal to $T_{U^0(0)} U^0(B(0,1))$, the tangent plane to the image of $U^0$ at $U^0(0)$. Let  $\delta>0$ and $r>0$ be small enough such that $C_{r}\cap U_0(A_\delta)$ is a simple curve. Then, for $\eps$ small enough, we easily see that the intersection of $u^\eps(A_\delta)$ and $C_r$ turn $l$ times around the cylinder, hence  $u^\eps(A_\delta)$ necessary intersect, which is a contradiction and proves the lemma.\hfill$\square$\\

{\bf Claim :
Let $u^\eps$ be a sequence of $C^2$-solutions of (\ref{princ}). We note $p$ the number of bubbles given by the decomposition \ref{de'}, this number splits into $k$ spheres and $l$ hemispheres, such that  $p=k+l$. Then we necessarily have that $l\geq 1$.}\\

{\it Proof of the claim :}\\

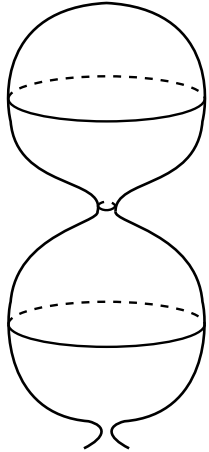
\begin{figure}[!h]
\centering
\input{explo1.tex}
\caption{Two bubbles joined by a neck}
\label{cou}
\end{figure}

We assume by contradiction that $l=0$. We show first that necessarily $k\leq 1$, that is to say there is no neck as in figure \ref{cou}.

We assume for contradiction that $k\geq 2$ and we consider the highest bubble, that is to say the one which  correspond to the smallest $\lambda_i^\eps$. Up to reorder, we assume that $i=1$. Thanks to (\ref{ortho'}), there is no bubble closed to this one, that is to say
\be
u^\eps_{1} \rightarrow \omega^1 \hbox{ in } C^2_{loc}(\R^2) 
\ee
where $u^\eps_{1}=u^\eps (a^\eps_1 + \lambda_1^\eps\, . \,)$.\\ 

We claim that the highest bubble is over another bubble. More precisely, there exists $i>1$ and $R_0>0$ such that  for all $R>0$ one get $B(a_1^\eps, R \lambda^\eps_1)\subset B(a_i^\eps, R_0 \lambda^\eps_i)$ for $\eps$ small enough.\\

Else this bubble would be isolated and will become tangent to $\partial \Omega$ at $0$. Indeed, there exists $z^0 \in \partial \D$ such that  for all $R>0$ and for all $z^\eps \in \partial B(a_1^\eps, R \lambda^\eps_1)$, there exists  a curve $\Gamma$ of $\overline{\D}$ joining $z^\eps$ to $z^0$ staying far from the other bubble. Hence thanks to the estimate (\ref{ii'}), for $\eps$ small enough, we see that the bubble $\omega^\eps_1$ would be almost tangent to $\partial\Omega^\eps$. This makes impossible the existence of an other isolated bubble, since it would be also almost tangent to $\partial\Omega^\eps$ at the same point, since we can take the same $z^0$ for the two bubbles. This would contradict the fact that the bubble are embedded and in the interior of $\Omega^\eps$. Hence the second bubble should be over the first one and so higher, which is a contradiction.\\

\begin{figure}[!h]
\centering
\input{explo2.tex}
\caption{$\tilde{\Sigma}^\eps$ range by  $u^\eps$ of a neighbourhood of the highest bubble}
\label{cou2}
\end{figure}
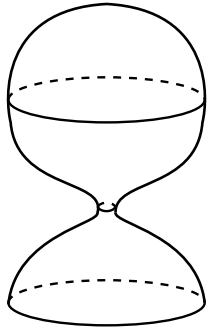

There exists $i_0$ and $R_0>0$ such that for all $R>0$ we get $B(a_1^\eps, R \lambda^\eps_1)\subset B(a_{i_0}^\eps, R_0 \lambda^\eps_{i_0})$. Then we choose the minimal  $\lambda_{i_0}^\eps$ satisfying this property. In this case we consider a neighbourhood  of $a^\eps_1$, $B(a_1^\eps, r \lambda_{i_0}^\eps)$,where $r>0$ is chosen such that this neighbourhood contains no other bubble. For $\eps$ small enough, the range of this neighbourhood by $u^\eps$, which will be noted $\tilde{\Sigma}^\eps$, seems like a sphere glued on a spherical cap, see figure \ref{cou2}.
 
Now we are in position to apply the Aleksandrov reflexion principle, as described in chapter VII of \cite{Hopf}. Let $P^\eps$ be the tangent plane to $\tilde{\Sigma}^\eps$ at $u^\eps(a^\eps_1)$ and $\nu^\eps$ the exterior normal at this point. We set $P^\eps_t =P^\eps + t \nu^\eps$, $\tilde{\Sigma}^{\eps,+}_t= \tilde{\Sigma}^\eps\cap \{ P^\eps + u \nu^\eps \hbox{ s.t. } u\geq t \}$ and $\tilde{\Sigma}^{\eps,-}_t$ the reflexion of  $\tilde{\Sigma}^{\eps,+}_t$ with respect to $P^\eps_t$. We consider the first negative time $t$ such that $\tilde{\Sigma}^{\eps,-}_t$ meets $\tilde{\Sigma}^{\eps}$. For $\eps$ small enough, the contact point of $\tilde{\Sigma}^{\eps,-}_t$ and $\tilde{\Sigma}^{\eps}$ can't belong to the boundary thanks to the presence of the neck. Moreover at the contact point the surfaces get the same orientation since  $\tilde{\Sigma}^{\eps,-}_t$ come from the interior of the bubble. Hence applying the Aleksendrov principle, we get $\tilde{\Sigma}^{\eps,-}_t = \tilde{\Sigma}^\eps \setminus \tilde{\Sigma}^{\eps,+}_t$, which is clearly a contradiction with the fact that the surface is embedded and with a boundary. Which proves that $k\leq 1$.\\

Now we have to exclude $k=0$ and $k=1$. If $k=0$, then $p=0$ and in that case the surface collapses. In fact his area will go to $0$, see the proof of theorem \ref{de'}. In that case, we rescale our space in order to get a new surface, denoted $\hat{\Sigma}^\eps$, whose area is equal to $1$. This imposes to the mean curvature of our new surface to goto $0$. Then our new sequence of surfaces  $\hat{\Sigma}^\eps$ goes to a minimal surface which bounds a plane curve. Indeed to insure the convergence it suffices to prove that $\vert \nabla \hat{u}^\eps \vert$ is uniformly bounded, where $\hat{u}^\eps$ is a conformal parametrization of $\hat{\Sigma^\eps}$. The regularity theory given in section \ref{sregul}, will give the convergence in $C^2(\overline{\D})$.\\

Let us assume by contradiction that $\ds \sup_\D \vert \nabla \hat{u}^\eps \vert \rightarrow +\infty$  when $\eps\rightarrow 0$. Then we set  $a^\eps_1\in \D$ et $\lambda_1^\eps$ such that 
\be \left\vert \nabla  \hat{u}^\eps  (a_{1}^\eps) \right\vert= \sup_{z\in \D}  \left\vert \nabla \hat{u}^\eps  (z) \right\vert 
\ee 
and
$$\left\vert \nabla \hat{u}^\eps  (a_{1}^\eps)\right\vert  = \frac{1}{\lambda_{1}^\eps} . $$
Then we set
$$\hat{u}_1^\eps (z)= \hat{u}^\eps(a_1^\eps +\lambda_1^\eps z )\hbox{ for all } z\in \Omega_1^\eps$$
where $\ds \Omega_1^\eps=\left\{ z\in R^2 \hbox{ s.t. } a_1^\eps +\lambda_1^\eps z\in \D\right\}$.\\

It is clear that $\vert \nabla \hat{u}^\eps_{1}\vert$ is bounded on every compact subset of $\overline{\Omega^\eps_1}$. Moreover, thanks to conformal invariance of our equation, $\hat{u}^\eps_{1}$ satisfies
$$ \Delta \hat{u}^\eps_{1} = o\left((\hat{u}^\eps_{1})_x \wedge( \hat{u}^\eps_{1})_y \right)$$
and
$$\langle (\hat{u}^\eps_{1})_x, (\hat{u}^\eps_{1})_y \rangle =\vert (\hat{u}^\eps_{1})_x\vert-\vert (\hat{u}^\eps_{1})_y\vert=0.$$
Hence, applying  standard elliptic theory, see section \ref{sregul} and \cite{GT}, we see that there exists a subsequence of $\hat{u}_1^\eps$ (still denoted  $\hat{u}_1^\eps)$ and $\beta^{1}\in C^2(\Omega_1)$ such that  
\be
\hat{u}^\eps_{1} \rightarrow \beta^{1} \hbox{ dans } C^2_{loc} (\overline{\Omega_{1}}) 
\ee
and
\be
\begin{split}
&\Delta \beta_1 =0,\\
&\langle (\beta_1)_x, (\beta_1)_y \rangle =\vert (\beta_1)_x\vert-\vert (\beta_1)_y\vert=0,
\end{split}
\ee
where $\displaystyle\Omega_1 =\lim_{\eps \rightarrow 0} \Omega^\eps_1$. Then there is two possibilities for $\Omega_1$, either it is the whole plane or it is a disk or an half-plane (which is conformaly equivalent). In this last case the boundary condition passes to the limit, that is to say, up to a rotation,
\be
\begin{split}
&\beta^{1}(\partial \Omega_{1}) \subset \R^2\times\{0\} \\
&\hbox { et }\\
& \left\langle  \beta^{1}_{x} \wedge \beta^{1}_{y}, N \right\rangle =0 \hbox{ on } \partial \Omega_{k} ,
 \end{split}
 \ee
where $N=(0,0,1)$. In that case $\beta_1$ can be extend by symmetry  in a $C^1$ function defined on the whole plane.\\
 
Moreover, thanks to the conformal invariance of $\Vert\nabla .\, \Vert_2$, we get 
$$\Vert \nabla \beta^{1} \Vert_2 \leq  2 \liminf_{\eps\rightarrow 0} \Vert \nabla \hat{u}^\eps_{1}\Vert_2= 2  \liminf_{\eps\rightarrow 0}\Vert \nabla \hat{u}^\eps\Vert_2=2.$$
Thanks to the Liouville theorem, we necessarily get $\nabla \beta_1 \equiv 0$, which is a contradiction with the fact that  $\vert \nabla \beta_{1}(0)\vert=1$. This proves that $\vert \nabla \hat{u}^\eps \vert$ is uniformly bounded and also the convergence of the sequence of surfaces $\hat{\Sigma}^\eps$.\\

With this convergence, the boundary condition passes to the limit, that is to say the minimal surface which is obtained meets the plane which contains its boundary orthogonally. But thanks to classical theory of minimal surfaces, see \cite{ColdingMinicozzi}, these surfaces should be flat, which contradicts the fact it must meet orthogonally the plane  which contains its boundary.\\

Finally, the last possibility is  $k=1$, that is to say there is only one bubble as in the following figure.
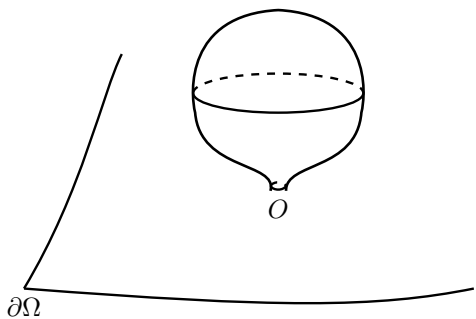
\begin{figure}[!h]
\centering
\input{explo3.tex}
\caption{a bubble meeting $\partial\Omega$}
\end{figure}  

But we can apply once more the Aleksandrov reflexion principle with respect to the tangent plane $\Sigma^\eps$ at the furthest point to $\partial \Omega^\eps$. Then the contact point between the surface and the reflected part is necessarily at the boundary, else the surface should be closed without boundary. Indeed if we have a contact in the interior, the local equality given by the reflexion principle would be global thanks to connexity, which is impossible since the upper part is simply connected and the lower is not. Hence the contact is done at the boundary $\partial\Sigma^\eps$, but  the tangent plane to  $\Sigma^\eps$ at the furthest point to  $\partial \Omega^\eps$ become parallel to the one of $\partial\Omega$ at $0$, which force the angle between $\partial\Omega^\eps$ and $\Sigma^\eps$ at the contact point to go zero when $\eps$ goes to $0$, which is a contradiction and achieves the proof of the claim.\hfill$\square$ 

\subsection{Proof of theorem \ref{Laurain10ter}}
Thanks to the previous section, $u^\eps(S^1)$ converges uniformly to a union of circles with radius $1$ centered at points $(c_i)$ of $T_0 \partial\Omega$. \\
\begin{figure}[!h]
\centering
\input{explo4.tex}
\caption{ $\omega^\eps$ which is bounded by $\partial \Sigma^\eps $ }
\end{figure}
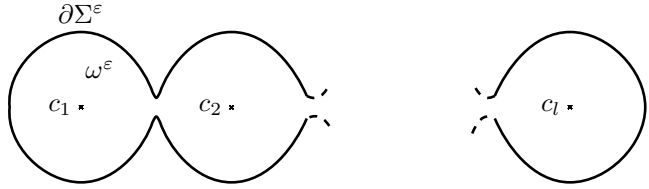  
Now we are in position to prove the theorem \ref{Laurain10ter}. In order to do it we apply the balancing formula (\ref{pbf3}) to the sequence $\Sigma^\eps$. This  gives
\beq
\label{bfth2}
2 \int_{\omega^\eps} \vec{N}^\eps dv -  \int_{\partial\Sigma^\eps} \vec{\nu}^\eps d\sigma=0
\eeq
where $\vec{\nu}^\eps$ is the conormal. \\

The fact that  $\Sigma^\eps$ and $\Omega^\eps$ meet orthogonally imposes $\vec{\nu}^\eps=\vec{N}^\eps$. We make a Taylor expansion of $\vec{N}^\eps$. Since $\partial \Omega^\eps$ is a graph above its tangent plane, we make the expansion in those coordinates. In fact, we are not going to do the expansion with respect to $0$ but with respect to $c^\eps\in\partial\Omega^\eps$ a point closed to $0$ which will be fixed later. 
\be
\vec{N}^\eps(z) = \vec{N}^\eps(c^\eps) +\eps d\vec{N}_{c^\eps}(z-c^\eps) + \eps^2 d^2\vec{N}_{c^\eps}(z-c^\eps)(z-c^\eps) +o(\eps^2).
\ee
At the first order, the left hand-side  (\ref{bfth2}) gives
\beq
\label{Xb}
\left(2 \vert \omega^\eps\vert  -\vert\partial\Sigma^\eps\vert\right)  \vec{N}(c^\eps).
\eeq
We can remark here that, thanks to theorem \ref{de'}, we have,  
$$ \vert \omega^\eps\vert\rightarrow l\pi,$$
and
$$\liminf_{\eps\rightarrow 0} \vert\partial\Sigma^\eps\vert \geq l2\pi .$$
In the first equality, we use the $L^\infty$-convergence while in the second, we use the $C^2_{loc}$-convergence. Then, thanks to (\ref{Xb}), we get 
$$\lim_{\eps\rightarrow 0} \vert\partial\Sigma^\eps\vert = l2\pi,$$
which proves that $\partial\Sigma^\eps$ converges to a union of $l$ circles as a current which justifies the fact that we will pass to the limit in the integral defined over this set.\\

In order to eliminate (\ref{Xb}) we project the left hand-side term of  (\ref{bfth2}) orthogonally to $\vec{N}^\eps({c^\eps})$, which gives to the second order
\beq
\label{last22}
\eps \pi^\eps\left(  2 \int_{\omega^\eps} d\vec{N}_{c^\eps}(z-c^\eps)\, dv -  \int_{\partial\Sigma^\eps}  d\vec{N}_{c^\eps}(z-c^\eps) \, d\sigma\right)
\eeq
where $\pi^\eps$ is the orthogonal projection parallel  to $\vec{N}^\eps({c^\eps})$.\\

Then we remark that there exists $c^\eps$ such that  (\ref{last22}) vanishes. Indeed
 \be 2 \int_{\omega^\eps} (z-c^\eps)\, dv -  \int_{\partial\Sigma^\eps}  (z-c^\eps) \, d\sigma
 \ee
is the weighted barycenter of $(\omega^\eps, 2)$ and $(\partial\Sigma^\eps, -1)$, then it suffices to choose $c^\eps$ as the corresponding barycenter to vanish (\ref{last22}).\\

Then it remains the order two terms, in which we pass to the limit after dividing them by $\eps^2$, which gives
\beq
\label{last23}
 \pi^0\left(  2 \int_{\omega^0} d^2\vec{N}_{0}(z-c^0) (z-c^0)\, dv -  \int_{\partial\Sigma^0}  d^2 \vec{N}_{c^0}(z-c^0)(z-c^0) \, d\sigma\right)=0
\eeq
where $\pi^0$ is the orthogonal projection  parallel to $\vec{N}(0)$, $c^0$, $\omega^0$ and $\partial\Sigma^0$ are respectively the limit of  $c^\eps$, $\omega^\eps$ and $\partial\Sigma^\eps$. As already remarked at the beginning of this section, $\omega^0$ and $\partial\Sigma^0$ are respectively an union of disks with radius $1$ and  union of circles with radius $1$ centered at some points $c_i$. Then we decompose the integral on this subset, which gives
\be
\begin{split}
&\sum_{i=1}^l \pi^0\left(  2 \int_{D(c_i,1)} d^2\vec{N}_{0}(z-c^0) (z-c^0)\, dv -  \int_{\partial D(c_i,1)}  d^2 \vec{N}_{c^0}(z-c^0)(z-c^0) \, d\sigma\right)=\\
&\sum_{i=1}^l \pi^0\left(  2 \int_{D(c_i,1)} d^2\vec{N}_{0}(z-c^i) (z-c^i)\, dv -  \int_{\partial D(c_i,1)}  d^2 \vec{N}_{c^0}(z-c^i)(z-c^i) \, d\sigma\right).
\end{split}
\ee
Here we use the fact that the integral vanishes if it contains an odd number of $z-c^i$ and the fact that $2\vert D(c_i,1)\vert = \vert\partial D(c_i,1)\vert$. Now we integrate, which traces $d^2 \vec{N}$ and gives
\be
l\left( 2 \int_\D \vert z\vert^2 \,dz - 1\right) \, \pi^0\left(\Delta \vec{N} (0)\right) =0.
\ee
But  the  general equation of a Gauss map of an immersion $X$ is given by
\beq
\label{fham}
\Delta \vec{N}=\vert \nabla \vec{N} \vert^2 \, \vec{N}  -2\vert \nabla X \vert^2\, \nabla H(u) ,
\eeq
hence we get
 \be
\frac{l}{2}\nabla H (0)=0,
\ee
which achieves the proof of the theorem.\hfill$\blacksquare$

\appendix

\section{Wente's inequality and application} 
\label{a5}
The aim of this section is to remind some Wente's inequalities, originally proved in  \cite{Wente69}. 

\begin{thm}
\label{wente1}
Let $\Omega$ be a bounded open set of $\R^2$ and $a,b\in H^1(\Omega)$. Let $u\in W^{1,1}_0(\Omega)$ be the solution of
$$\Delta u =  a_x b_y -a_yb_x \hbox{ on } \Omega ,$$
then 
$$\Vert u\Vert_{\infty} \leq \frac{1}{2\pi} \Vert \nabla a\Vert_2 \Vert \nabla b\Vert_2 ,$$ 
and 
$$\Vert \nabla u\Vert_{2} \leq \sqrt{ \frac{3}{16\pi}} \Vert \nabla a\Vert_2 \Vert \nabla b\Vert_2 .$$
Moreover the constant are optimal.
\end{thm}
Which is remarkable here is that the constant is independent of  $\Omega$. We will find the proof in  \cite{Topping97} and \cite{Ge98}, see also \cite{Wente80} and \cite{Baraket96}. These inequalities have been extend to function defined on surfaces. In particular, we have the following theorem. 
\begin{thm}
\label{topp} Let $\Sigma$ a compact Riemannian surface without boundary and $v\in H^1(\Sigma, \R^2)$. Then if $u\in W^{1,1}(\Sigma,\R)$ be the solution of
$$\Delta u = det(\nabla v) \hbox{ on } \Sigma ,$$
then 
$$osc(u)+\Vert \nabla u \Vert_2 \leq \left( \frac{1}{4\pi}+  \sqrt{ \frac{3}{128\pi}} \right)  \Vert \nabla v\Vert_2^2 ,$$
where $\ds osc(u)= \sup_{x,y\in \Sigma}\vert u(x)-u(y)\vert .$ 
\end{thm}
Then,  assuming that $u\in H^1$, we extend such an equality to $\Omega=\R^2$.
\begin{cor}
\label{wente2}
Let $v\in H^1(\R^2,\R^3)$ and  $u\in H^{1}(\R^2,\R^3)$ be a solution of
$$\Delta u = -2  v_x \wedge v_y \hbox{ on } \R^2$$
then
$$osc(u) +\Vert \nabla u \Vert_2 \leq  \left( \frac{1}{\pi}+  \sqrt{ \frac{3}{ 8\pi}} \right)\Vert \nabla v\Vert_2^2 .$$ 
\end{cor}
Here the constant is {\it a priori} not optimal.\\

{\it Proof of corollary \ref{wente2} :}\\

Let $\pi$ the standard stereographic projection from $S^2$ to $\R^2$. Thanks to the conformal invariance of the equation, $u \circ \pi^{-1}$ and $v\circ \pi^{-1}$ satisfies the hypothesis of  theorem \ref{topp} when $\Sigma=S^2$, hence we get that 
$$osc(u^1)+ \Vert \nabla u^1\Vert_2 \leq \left( \frac{1}{2\pi} + \sqrt{\frac{3}{32\pi}}\right) \left(\Vert \nabla v^2 \Vert_2^2 + \Vert \nabla v^3 \Vert_2^2\right) , $$
$$osc(u^2)+ \Vert \nabla u^3\Vert_2 \leq \left( \frac{1}{2\pi} + \sqrt{\frac{3}{32\pi}}\right) \left(\Vert \nabla v^1 \Vert_2^2 + \Vert \nabla v^3 \Vert_2^2\right) $$
and 
$$osc(u^3)+ \Vert \nabla u^3\Vert_2 \leq \left( \frac{1}{2\pi} + \sqrt{\frac{3}{32\pi}}\right) \left(\Vert \nabla v^1 \Vert_2^2 + \Vert \nabla v^2 \Vert_2^2\right) .$$

Then summing this inequalities, we get the desired inequality.\hfill$\square$\\

To conclude this section we remind a useful  Wente's type inequality, see \cite{BC1} for example.
\begin{lemma}
\label{Xlem} 
Let  $u\in H^1(\D)\cap L^\infty(\D)$ and $v\in H^1_0(\D)$, then there exists $C$, independent of $u$ and $v$, such that 
$$\left\vert \int_{\Omega} \langle u, v_x \wedge v_y \rangle \right\vert\leq C \Vert \nabla u\Vert_2  \Vert \nabla v\Vert_2^2 .$$
\end{lemma}

\addcontentsline{toc}{chapter}{Bibliographie}
\bibliographystyle{plain}
\bibliography{cmcbiblio}

\end{document}

%% file: figbf.tex
\begin{tikzpicture}[scale=0.5, arrow/.style={postaction={decorate}, decoration={markings, mark= at position 1  with {\arrow[line width=1pt]{>}}}}]] 

\draw[dashed, line width=1pt] (-5.5,-0.2) -- (-4.7,0);

\draw[line width=1pt] (-4.7,0).. controls (-2,1) and( -0,1) .. (2,0) .. controls (4,-1) and (5,-1) .. (7,0);

\draw[dashed, line width=1pt] (7,0) -- (8.7,0.7);

\draw[ line width=1pt, arrow] (7,0) -- (8,0.4);

\draw (7.3,0.1) -- (7.5,.35);
\draw (7.5,.35) -- (7.1,0.2);

\draw (7.5,0) node[below]{$\vec{\eta}$} ;

\draw[dashed, line width=1pt] (-5.4,2.7) -- (-4.7,3);

\draw[line width=1pt] (-4.7,3).. controls (-2,4) and( -1,3)  .. (0.3,5).. controls (1,6) and (2,7) .. (3,6.5) .. controls (4,6) and (5,5) .. (7,4);

\draw[dashed, line width=1pt] (7,4) --(8.7,3.2);

\draw[line width=1pt] (7,4).. controls (6,3.5) and( 6,0) .. (7,-0);
\draw[dashed,line width=1pt] (7,4).. controls (8,3.8) and( 8,0.5) .. (7,0);

\draw[line width=1pt] (-4.7,0).. controls (-5.5,0.5) and( -5.5,2.5) .. (-4.7,3);
\draw[dashed,line width=1pt] (-4.7,0).. controls (-4,0) and( -4,3) .. (-4.7,3);

\draw[line width=1pt] (1.15,6).. controls (2,5.5) and( 3.2,5.5) .. (3.8,6);
\draw[dashed,line width =1pt] (1.15,6).. controls (2.1,6.2) and( 3.1,6.2) .. (3.8,6);

\draw (0,2) node[above]{$S$} ;

\fill[color=gray!20] (7,4).. controls (6,3.5) and( 6,0) .. (7,0) .. controls (8,0.5) and( 8,3.8)..(7,4);
\fill[color=gray!20] (-4.7,0).. controls (-5.5,0.5) and( -5.5,2.5) .. (-4.7,3) .. controls (-4,3) and( -4,0) ..(-4.7,0);

\draw (-4.7,1.5) node[above]{$\Sigma$} ;
\draw (7,2) node[above]{$\Sigma$} ; 

\draw[ line width=1pt, arrow] (7,2) -- (8,2);

\draw (7.2,2) -- (7.2,1.8);
\draw (7.2,1.8) -- (7,1.8);

\draw (8,2) node[right]{$\vec{\nu}$} ;

\draw (-5.8,2) node {$\partial S$} ; 
\draw (4.9,2) node[right]{$\partial S$} ;

\draw (0,4.8)--(2,7);
\draw[ line width=1pt, arrow] (0.9,5.8) -- (0.1,6.6);

\draw (0,6.7) node[above]{$\vec{N}$} ;

\draw (1.05,5.9) -- (0.87,6.1);
\draw (0.87,6.1) -- (0.75,5.95);

\end{tikzpicture}

%% file: explo1.tex
\begin{tikzpicture}[scale=0.3] 

\draw[line width=1pt] (0.4,5) .. controls (0.35, 6 .5) and (4,6) .. (4.25,9) .. controls (4.5,10) and (4.75,14) .. (0,14.25 ) .. controls (-4.75,14) and (-4.5,10) .. (-4.25,9) .. controls  (-4, 6) and (0, 6.5) .. (-0.4,5).. controls  (-0.3, 4.5) and (-4,4) .. (-4.25,1) .. controls (-4.5,0) and (-4.5,-4) .. (-0.8,-4.3 ) ..
controls (0,-4.5) and (0,-5) .. (-1,-5.5 ) ;

\draw[line width=1pt]  (1,-5.5 )..controls (0,-5) and (0,-4.5) .. (0.8,-4.3 ).. controls (4.5,-4) and (4.5,0) .. (4.25,1) .. controls  (4, 4) and (0.3, 4.5) .. (0.4,5)  ;

\draw[line width=0.9pt] (-4.35,10) arc (-180:0:4.35cm and 1cm);
\draw[dashed,line width=0.9pt] (-4.35,10) arc (180:0:4.35cm and 1cm);

\draw[line width=0.9pt] (-4.35,0) arc (-180:0:4.35cm and 1cm);
\draw[dashed,line width=0.9pt] (-4.35,0) arc (180:0:4.35cm and 1cm);

\draw[line width=0.9pt]  (-0.4,5.25) ..controls (-0.25,5) and (0.25, 5) .. (0.4,5.25);
\draw[dashed,line width=0.9pt]  (-0.4,5.25) ..controls (-0.25,5.5) and (0.25,5.5) .. (0.4,5.25);

\end{tikzpicture}

%% file: explo2.tex
\begin{tikzpicture}[scale=0.3] 

\draw[line width=1pt] (0.4,5) .. controls (0.35, 6 .5) and (4,6) .. (4.25,9) .. controls (4.5,10) and (4.75,14) .. (0,14.25 ) .. controls (-4.75,14) and (-4.5,10) .. (-4.25,9) .. controls  (-4, 6) and (0, 6.5) .. (-0.4,5).. controls  (-0.3, 4.5) and (-4,4) .. (-4.35,0.95);

 \draw [line width=1pt] (4.35,1) .. controls  (4, 4) and (0.3, 4.5) .. (0.4,5)  ;

\draw[line width=0.9pt] (-4.35,10) arc (-180:0:4.35cm and 1cm);
\draw[dashed,line width=0.9pt] (-4.35,10) arc (180:0:4.35cm and 1cm);

\draw[line width=0.9pt] (-4.35,1) arc (-180:0:4.35cm and 1cm);
\draw[dashed,line width=0.9pt] (-4.35,1) arc (180:0:4.35cm and 1cm);

\draw[line width=0.9pt]  (-0.4,5.25) ..controls (-0.25,5) and (0.25, 5) .. (0.4,5.25);
\draw[dashed,line width=0.9pt]  (-0.4,5.25) ..controls (-0.25,5.5) and (0.25,5.5) .. (0.4,5.25);

\end{tikzpicture}

%% file: explo3.tex
\begin{tikzpicture}[scale=0.26] 

\draw[line width=1pt] (0.4,5) .. controls (0.35, 6 .5) and (4,6) .. (4.25,9) .. controls (4.5,10) and (4.75,14) .. (0,14.25 ) .. controls (-4.75,14) and (-4.5,10) .. (-4.25,9) .. controls  (-4, 6) and (0, 6.5) .. (-0.4,5);

\draw[line width=0.9pt] (-4.35,10) arc (-180:0:4.35cm and 1cm);
\draw[dashed,line width=0.9pt] (-4.35,10) arc (180:0:4.35cm and 1cm);

\draw[line width=0.9pt]  (-0.4,5.25) ..controls (-0.25,5) and (0.25, 5) .. (0.4,5.25);
\draw[dashed,line width=0.9pt]  (-0.4,5.25) ..controls (-0.25,5.5) and (0.25,5.5) .. (0.4,5.25);

\draw[line width=0.9pt]  (10,0) .. controls (5,-1) and (0,-1) .. (-13,0); 
\draw[line width=0.9pt]  (-13,0) .. controls (-10,5) and (-9,10) .. (-8,12); 

\draw (-13,0) node[below]{$\partial\Omega$} ;

\draw (0,5) node[below]{$O$} ;

\end{tikzpicture}

%% file: explo4.tex
\begin{tikzpicture}[scale=0.5] 

\draw[dashed, line width=1pt] (0.6,4.5) .. controls (0.4,4.8) and (0.2,4.8) ..(-0,4.7);

\draw[line width=1pt] (-0,4.7).. controls (-0.3,4) and (-1,3)..(-2,3).. controls (-3,3) and (-3.7,4).. (-3.9,4.6).. controls (-4,4.8) and (-4,4.8)..(-4.1,4.6).. controls (-4.3,4) and (-5,3)..(-6,3).. controls (-7,3) and (-8,4)..(-7.9,5).. controls (-8,6) and (-7,7).. (-6,7) .. controls (-5,7) and  (-4.3,6).. (-4.1,5.4)..controls(-4,5.2) and (-4,5.2) .. (-3.9, 5.4) .. controls (-3.7,6) and (-3,7).. (-2,7) .. controls (-1,7) and  (-0.3,6) .. (0, 5.3);

\draw[dashed, line width=1pt] (0,5.3) .. controls (0.2,5.2) and (0.4,5.2) ..(0.6,5.6);

\draw[dashed, line width=1pt] (4.4,4.3) .. controls (4.6,4.8) and (4.8,4.8) ..(5,4.7);

\draw[line width=1pt] (5,4.7).. controls (5.3,4) and (6,3)..(7,3).. controls (8,3) and (9,4).. (9,5).. controls (9,6) and (8,7).. (7,7) .. controls (6,7) and  (5.3,6) .. (5, 5.3);

\draw[dashed, line width=1pt] (5,5.3) .. controls (4.8,5.2) and (4.6,5.2) ..(4.4,5.7);

\draw (-5.5,6) node {$\omega^\eps$} ;
\draw plot[mark=x,mark size=0.7 mm] (-6,5) node[left] {$c_1$} ;
\draw plot[mark=x,mark size=0.7 mm] (-2,5) node[left] {$c_2$} ;
\draw plot[mark=x,mark size=0.7 mm] (7,5) node[left] {$c_l$} ;

\draw (-6,7) node[above] {$\partial\Sigma^\eps$} ;

\end{tikzpicture}